\documentclass[10pt,twoside]{amsart}
\usepackage{amssymb,amsbsy,amsmath,amsfonts,amssymb,amscd}
\usepackage{latexsym,euscript,exscale}%eucal,epic,eepic,epsfig
\usepackage{times}

\newcommand{\emptyprop}{q}

\newcommand \map[1]{{\newcommand{\tmpprop}{#1q}  \if\tmpprop\emptyprop \to\else \xrightarrow{{\phantom{i}{#1}\phantom{i}}}\fi}} 
\newcommand \maxim{\mathfrak m}
\newcommand \nat{\mathbb N}

\newcommand \pol[2]{#1[#2]}
\newcommand \pow[2]{#1[[#2]]}

\newcommand \tor[4]{\operatorname{Tor}^{#1}_{#2}(#3,#4)}
\newcommand \op\operatorname

\newcommand \ch{characteristic}
\newcommand \homo{homomorphism}
\newcommand \CM{Coh\-en-Mac\-au\-lay}
\renewcommand\iff{if, and only if,}
\newcommand \DVR{discrete valuation ring}

\newtheorem{theorem}{Theorem}[section]

\newtheorem{e-proposition}[theorem]{Proposition}

\newtheorem{e-definition}[theorem]{Definition\rm}

\title {Mixed characteristic homological theorems in low degrees}
\author{Hans Schoutens}
\thanks{Partially supported by a  grant from the National Science Foundation.}
\address{Department of Mathematics\\
Ohio State University\\
Columbus, OH 43210 (USA)}
\email{schoutens@math.ohio-state.edu}

\begin{document}

\begin{abstract}   
Let  $R$  be a locally finitely generated algebra over a \DVR\ $V$ of mixed \ch. For any of the homological properties, the Direct Summand Theorem, the Monomial Theorem, the Improved New Intersection Theorem, the Vanishing of Maps of Tors and the Hochster-Roberts Theorem, we show that it holds for $R$ and possibly some other data defined over $R$, provided the residual \ch\ of $V$ is sufficiently large in terms of the complexity of  the data, where the complexity is primarily given in terms of the degrees of the polynomials over $V$ that define the data, but possibly  also by some additional invariants.
\end{abstract}

\maketitle

\newcommand\preprintnote {preprint on \myhomepage}
\newcommand\myhomepage{http://www.math.ohio-state.edu/\-\~{}schoutens}

%LOCAL DEFINITIONS

\newcommand \zet{\mathbb Z}
\newcommand \ul[1]{\seq{#1}\infty}

\newcommand \seq[2]{#1\mathstrut_{#2}}
\newcommand \sr{approximation}
\newcommand \uleq[1]{\ulseq{#1}\infty{eq}}
\newcommand \ulseq[3]{#1\mathstrut^{\text{#3}}_{#2}}
\newcommand \ulmix[1]{\ulseq{#1}\infty{mix}}
\newcommand \los{\L os' Theorem}
\newcommand \br{\mathfrak O}

\newcommand  \q{pseudo}
\newcommand  \qdim{\q-dimension}
\newcommand  \qCM{\q-\CM}
\newcommand  \qreg{\q-regular}

\section{The results}

Let $V$ be a mixed \ch\ \DVR\ with uniformizing parameter $\pi$ and residue field $\kappa$ of \ch\ $p$. We say that $R$ is a \emph{local $V$-affine algebra} of \emph{$V$-complexity} at most $c$, if it is of the form $(\pol VX/I)_\maxim$, with $X$ a tuple of at most $c$ variables, $I$ and $\maxim$ ideals generated by polynomials of degree at most $c$, and $\maxim$ a prime ideal containing $I$ and   $\pi$. Similarly, we say that an element $a$ in $R$ (respectively, a tuple $\mathbf x$ in $R$; a matrix $\Gamma$ defined over $R$; an ideal  $I$ in $R$; a finitely generated $R$-module $M$; or, an $R$-algebra  $S$) has $V$-complexity at most $c$, if $R$ has $V$-complexity at most $c$ and $a$ is the image in $R$ of a fraction $f/g$ with $f$ and $g$ polynomials of degree at most $c$ and $g\notin\maxim$ (respectively, the length of $\mathbf x$ is at most $c$ and each of its entries has $V$-complexity at most $c$; the dimensions of $\Gamma$ are at most $c$ and each of its entries has $V$-complexity at most $c$; the ideal $I$ is generated by elements of $V$-complexity at most $c$; the module $M$ can be realized as the cokernel of a matrix of $V$-complexity at most $c$;  and,  the $R$-algebra $S$ has $V$-complexity  at most $c$).

\begin{theorem}[Asymptotic Homological $\mathcal P$-Theorem]\label{T:main}
Let  $\mathcal P$ be one of the homological properties listed below. For each $c\in\nat$, there exists a bound $c'\in\nat$, such that if  $V$ is a  mixed \ch\ \DVR, $R$ a local $V$-affine algebra, and  $\varpi$ some other data defined over $R$, all of $V$-complexity at most $c$ (and possibly with some additional constraints in terms of $c$ indicated below),  and if the residual \ch\  of $V$ is at least $c'$, then  property $\mathcal P$ holds for $\varpi$.

\begin{description}
\item[Direct Summand Theorem.] Given a module-finite ring extension $R\subset S$, if $R$ is regular, then $R\subset S$ splits as an $R$-module morphism.
\item[Monomial Theorem.] Given at most $c$ monomials $Y^{\mu_i}$ in at most $c$  variables $Y$ and given a system of parameters $\mathbf x$ of $R$, such that $\mathbf xR\cap V$ has  $V$-adic valuation at most $c$, if $Y^{\mu_0}$ does not belong to the ideal in $\pol\zet Y$ generated by the remaining monomials $Y^{\mu_i}$, then  $\mathbf x^{\mu_0}$ does not belong to the ideal in $R$ generated by the remaining $\mathbf x^{\mu_i}$. 
\item[Improved New Intersection Theorem.] Given a finite free complex 
	\begin{equation}
	 0\to R^{a_s}\map{\Gamma_s} R^{a_{s-1}}\map{\Gamma_{s-1}}\dots \map{\Gamma_2} R^{a_1} \map{\Gamma_1} R^{a_0}\to 0\tag{$F_\bullet$}
	\end{equation}
over $R$ with $s,a_i\leq c$ and a minimal generator $\tau$ of $H_0(F_\bullet)$ generating a module of length at most $c$, if each $R/I_{r_i}(\Gamma_i)$ has  dimension at most $d-i$ and parameter degree\footnote{The \emph{parameter degree} of a Noetherian local ring $S$ is defined as the minimal possible length of a residue ring $S/\mathbf xS$, where $\mathbf x$ runs over all systems of parameters of $S$ (note that homological multiplicity is an upper bound for parameter degree by \cite[\S4]{SchABCM}). } at most $c$, where  
	\begin{equation*}
	r_i:= \sum_{j=i}^s (-1)^{j-i} a_j,
	\end{equation*}
and $d$ is the dimension of $R$, then $F_\bullet$ has length at least $d$. Here we write $I_n(\Gamma)$ for the  ideal generated by all $n\times n$-minors of a matrix $\Gamma$.
\item[Vanishing for Maps of Tors.] Given $V$-algebra \homo{s} $R\to S\to T$  and a finitely generated $R$-module $M$, if $R$ and $T$ are regular and if $R\to S$ is integral and injective, then the natural map
	\begin{equation*}
	\tor RnSM \to \tor RnTM
	\end{equation*}
is zero.
\item[Hochster-Roberts Theorem.] Given a cyclically pure\footnote{A \homo\ $R\to S$ is \emph{cyclically pure} if $I=IS\cap R$, for every ideal $I$ in $R$.}  \homo\ $R\to S$ of $V$-algebras, if $S$ is regular, then $R$ is \CM.
\end{description}
\end{theorem}

\section{The method}

If $V$ is equi\ch, then each of these homological  properties  holds unconditionally, that is to say, without any bound on the complexity (\cite{BH,HoDS,Str90}). We will use the Ax-Kochen-Ershov Principle to deduce Theorem~\ref{T:main} from this. Let me sketch the idea before I give more details. After a faithfully flat extension, we may assume that $V$ is moreover complete. Towards a contradiction, suppose for some $c$, no such bound exists. This means that for each $p$, we can find a complete \DVR\ $\seq Vp$ of \ch\ zero and  residual \ch\ $p$, and  some data $\seq\varpi p$ of $\seq Vp$-complexity at most $c$ for which $\mathcal P$ fails. Let $\seq\kappa p$ be the residue field of $\seq Vp$. Define $\ulseq Vp{eq}:=\pow{\seq\kappa p}t$, for $t$ a single variable. Using the Ax-Kochen-Ershov Principle, we can construct for each $p$, similar data $\ulseq\varpi p{eq}$ defined over the \DVR{s} $\ulseq Vp{eq}$, so that for infinitely many $p$, property $\mathcal P$ does not hold for $\ulseq\varpi p{eq}$, leading to the desired contradiction. 

I will now explain this in more detail. The relation between the \DVR{s} $\seq Vp$ and $\ulseq Vp{eq}$ is given by the following result due to Ax-Kochen \cite{AK} and Ershov \cite{Ers65,Ers66}.

\begin{theorem}[Ax-Kochen-Ershov]
For a fixed choice of a non-principal ultrafilter on the set of prime numbers, the ultraproduct  of all $\seq Vp$ is isomorphic to the ultraproduct of all $\ulseq Vp{eq}$.
\end{theorem}

For a quick review on ultraproducts, including \los, see \cite[\S2]{SchNSTC}; for a more detailed treatment, see \cite{Hod}. Fix a non-principal ultrafilter on the set of prime numbers. Identify both ultraproducts via a fixed isomorphism and denote the common object by $\br$. By \los, $\br$ is an equi\ch\ zero Henselian (non-discrete, non-Noetherian) valuation ring with maximal ideal generated by a single element $\pi$. Fix a tuple of variables $X$. It is no longer true that the ultraproduct $\ulmix A$ of the $\pol{\seq Vp}X$ is isomorphic to the ultraproduct $\uleq A$ of the $\pol{\ulseq Vp{eq}}X$. Nonetheless, both ultraproducts contain $\pol \br X$ as a subring. More precisely, if $\seq fp\in\pol{\seq Vp}X$ have degree at most $c$, for some $c$ independent from $p$, then their ultraproduct $\ul f$ in $\ulmix A$ is an element of the subring $\pol \br X$, and every element in $\pol \br X$ is realized in this manner. In particular, $\ul f$ can also be viewed as an element in $\uleq A$, that is to say, as the ultraproduct of elements $\ulseq fp{eq}\in\pol{\ulseq Vp{eq}}X$. In this way, we can associate to a sequence of elements $\seq fp$ of uniformly bounded $\seq Vp$-complexity, a sequence of elements $\ulseq fp{eq}$ of uniformly bounded $\ulseq Vp{eq}$-complexity. Although this assignment is not unique, any two choices will be the same almost everywhere (in the sense of the ultrafilter). Similarly, we can associate to a sequence of local $\seq Vp$-affine algebras $\seq Rp$ of uniformly bounded complexity (or any other object defined in finite terms over $\seq Vp$), a sequence of local $\ulseq Vp{eq}$-affine algebras $\ulseq Rp{eq}$; the latter are called an equi\ch\ \emph{\sr} of the former.

Let $\ulmix R$ and $\uleq R$ be the respective ultraproducts of $\seq Rp$ and $\ulseq Rp{eq}$. These rings have a common local subring $(\Re,\maxim)$, consisting precisely of ultraproducts of elements of uniformly bounded complexity. Using for instance the result in \cite{Asch} regarding uniform bounds on the complexity of modules of syzygies, one shows that both extensions $\Re\to \ulmix R$ and $\Re\to \uleq R$  are faithfully flat. Moreover,  using results from \cite{SvdD,SchBC,SchBounds}, every finitely generated prime ideal of $\Re$ remains prime when extended to either $\ulmix R$ or $\uleq R$.  It follows that almost all $\seq Rp$ are domains \iff\ $\Re$ is a domain \iff\ almost all $\ulseq Rp{eq}$ are domains.

 The idea is to view $\Re$ as an equi\ch\ zero version of the $\seq Rp$ (or, for that matter, of the $\ulseq Rp{eq}$), so that we are lead to prove an analogue of the homological property $\mathcal P$ for $\Re$ (and whatever other data required, arising in a similar fashion from data of uniformly bounded $\seq Vp$-complexity). However, in carrying out this project, we are faced with a serious obstruction: $\Re$  is in general not Noetherian. This prompts for a non-Noetherian version of the local algebra required for discussing homological properties. To this end, we define the \emph{\qdim} of $\Re$ to be the smallest length of a tuple generating an $\maxim$-primary ideal (note that the Krull dimension is infinite and hence of no use). We say that $\Re$ is \emph{\qreg} if its \qdim\ equals its embedding dimension (=the minimal number of generators of $\maxim$), and \emph\qCM, if its \qdim\ is equal to its depth (in the sense of \cite{HoGrade}). To derive for instance the asymptotic Hochster-Roberts Theorem, we can now use the fact that almost all $\seq Rp$ are regular (respectively, \CM) \iff\ $\Re$ is \qreg\ (respectively, \qCM)  \iff\ almost all $\ulseq Rp{eq}$ are regular (respectively, \CM).

The main tool in establishing a variant of each $\mathcal P$ over $\br$ is via an $\br$-analogue of a big \CM\ algebra. Hochster has demonstrated (see for instance \cite{HoHT,HoDS}) how efficiently  big \CM\ modules can be used to prove homological theorems. More recently, Hochster and Huneke have given various strengthenings and generalizations using big \CM\ \emph{algebras}.  Big \CM\ algebras in equi\ch\ zero are obtained by reduction from their existence in \ch\ $p$ via absolute integral closures (\cite{HHbigCM,HHbigCM2}). In \cite{SchBCM}, I gave an alternative construction for local $\mathbb C$-affine domains, using ultraproducts, and it is this approach we will adopt here. Namely, for $\Re$ a local $\br$-affine domain, let $\mathcal B(\Re)$ be the ultraproduct of the absolute integral closures $(\ulseq Rp{eq})^+$.

\begin{theorem}[Big \CM\ Algebra]
Let $(\Re,\maxim)$ be a local $\br$-affine domain.  Every tuple  in $\Re$ of length equal to the \qdim\ of $\Re$ and generating an $\maxim$-primary ideal, is $\mathcal B(\Re)$-regular.
\end{theorem}
\begin{proof}
Let $\mathbf x$ be a tuple of length equal to the \qdim\  $d$ of $\Re$ so that $\mathbf x\Re$ is $\maxim$-primary. Choose $d$-tuples $\ulseq{\mathbf x}p{eq}$ in $\ulseq Rp{eq}$ whose ultraproduct is $\mathbf x$. One can show that almost all $\ulseq Rp{eq}$ have dimension $d$. By \los, almost all $\ulseq{\mathbf x}p{eq}\ulseq Rp{eq}$ are primary to the maximal ideal. Hence almost all $\ulseq{\mathbf x}p{eq}$ are systems of parameters, whence $(\ulseq Rp{eq})^+$-regular by \cite{HHbigCM}. By another application of \los, $\mathbf x$ is $\mathcal B(\Re)$-regular.
\end{proof}

Details can be found in the forthcoming \cite{SchMixBCM}.

%%%%%%%%%%%%%%%%
%%%  Bibliography  %%%
%%%%%%%%%%%%%%%%%%%%%%%

\providecommand{\bysame}{\leavevmode\hbox to3em{\hrulefill}\thinspace}

\end{document}